\documentclass{birkjour}
\usepackage{amsmath}
\usepackage{amssymb}
\begin{document}
\title{Notes on  Chern's \\  Affine Bernstein Conjecture}
\author{An-Min Li}
\address{Yangtze Center of Mathematics, Sichuan University,
        Chengdu,  P.R. China}
\email{math$\_$li\@yahoo.com.cn}
\author{Ruiwei Xu}
\address{Department of Mathematics, Henan Normal University,
        Xinxiang,  P.R. China}
        \email{xuruiwei@yahoo.com.cn}
\author{Udo Simon}
\address{Institut f\"ur  Mathematik, TU Berlin,
        Germany}
        \email{simon@math.tu-berlin.de}
\author{Fang Jia }
\address{Department of Mathematics, Sichuan University,
        Chengdu,  P.R. China}
\email{galaxyly2000@yahoo.com.cn}

\begin{abstract} There were two famous conjectures on
{\it complete affine maximal surfaces}, one due to E. Calabi, the
other to S.S. Chern.  Both were solved with different methods about
one decade ago by studying the associated Euler-Lagrange equation.
Here we survey  two proofs of Chern's conjecture in our recent
monograph \cite{L-X-S-J}, in particular we add some  details of the
proofs of auxiliary material that were omitted in \cite{L-X-S-J}. We
describe the related background  in our Introduction. Our survey is
suitable as a report about recent developments and techniques in the
study of certain Monge-Amp\`{e}re equations.
\end{abstract}
\subjclass{53\,A\,15, 35\,J\,60, 58\,J\,60.} \keywords{Chern's
Affine Bernstein Conjecture, affine maximal surfaces,
Monge-Amp\`{e}re equations, affine blow-up analysis} \maketitle
\section{ Introduction}
Many geometric problems in analytic formulation lead to important
classes of PDEs. Naturally, geometric methods play a crucial role in
their investigation. Typical examples are the {\it Euclidean
Minkowski Problem}, the Theorem of J\"orgens-Calabi-Pogorelov on the
{\it global uniqueness of improper affine spheres}, or the two
famous {\it Affine Bernstein Conjectures} of Calabi and Chern, resp.\\
Wide use of geometric methods in studying PDE's of affine
hypersurface theory   was initiated by E. Calabi and  continued by
A.V. Pogorelov, S.Y. Cheng - S.T. Yau, A.-M. Li, and, during the
last decade,  e.g. by N.S. Trudinger - X.J. Wang,  A.-M. Li's
school, and other authors.  An effective application  of such
geometric methods has played  a crucial role in these studies. \\
The contributions of E. Calabi and  S.Y. Cheng - S.T. Yau had a
particularly deep influence on the development of this subject.
According to the foreword in \cite{Ch-Y-86} this paper originated
from discussions with E. Calabi and L. Nirenberg, and from results
of both on the same topic; for further historical details and
references we refer to
\cite{Ca}, \cite{CA}, \cite{POG-72}, \cite{POG-75}
and the  monographs \cite{L-S-Z},  \cite{L-X-S-J}.\\
In problems involving PDE's of Monge-Amp\`{e}re type it is often the
case that the unknown solution is a convex function defining locally
a nonparametric hypersurface for which it is possible to choose a
suitable relative normalization and investigate the induced
geometry. We refer to this process as {\it geometric modelling},
described and emphasized in the introduction of the monograph \cite{L-X-S-J}.\\
Our emphasis is different here. Even when there is some parallelity
in tools, that is Lemmas and Propositions, used in \cite{L-X-S-J}
and also here, here our emphasis is on the following:\\
In recent  years, A.-M. Li and his school extended the study of
certain types of Monge-Amp\`{e}re equations and developed  a
framework, in particular  for studying such types of 4-th order
PDE's, including the affine maximal hypersurface equation and the
Abreu equation (cf. \cite{L-J-1}, \cite{L-J-2}, \cite{L-J-3},
\cite{L-J-4}, \cite{L-X-1},\cite{L-X-2}, \cite{C-L-S-2}, \cite{C-L-S-4}).\\
We call this the {\em real affine technique}. The whole package
includes:
\begin{itemize}
\item the derivation of
differential inequalities for certain functions of \\geometric
importance, related to the given problem;  here that is  the
differential inequality (4.3) for the function $\Phi$;
\item convergence
theorems,
\item Bernstein properties and
\item an affine blow-up analysis.
\end{itemize}
This package is very useful for studying  such types of PDEs; the
technique was extended to complex manifolds; it also plays an
important role in the study of extremal metrics on toric surfaces;
see \cite{C-L-S-1} and \cite{C-L-S-3}.\\
It is the central  aim of this paper to survey this {\em real affine
technique} and to sketch  two proofs of Chern's conjecture in our
recent monograph \cite{L-X-S-J}, in particular to outline details
of the proof of Proposition 5.6.13, which were  omitted in
\cite{L-X-S-J}; see section 3 below.

{\bf Chern's Conjecture.} {\it $\,$ Let
$\;x_{n+1}=f(x_1,x_2,...,x_n)$ be a smooth, strictly convex function
defined for all $(x_1,x_2,...,x_n)\in \mathbb{R}^n$. If the graph
hypersurface
$$M=\{(x,f(x)\;|\; x\in \mathbb{R}^n\}$$ is affine maximal
then $M$ must be an elliptic paraboloid.}

The two-dimensional conjecture of Chern was solved by N. Trudinger
and X. Wang in \cite{T-W-1}; they combined {\it tools from  the theory of
convex bodies}   and from the {\it Caffarelli-Guti\'{e}rrez theory}. In
\cite{L-J-3} one can find a  different proof, also using tools from the theory
of convex bodies  and the Caffarelli-Guti\'{e}rrez theory.
In \cite{L-J-1} and \cite{L-X-S-J} the authors gave a completely
different proof of Chern's conjecture, purely using tools from analysis.\\
Concerning this purely  analytic  proof of Chern's conjecture, our proof
here  also meets demands of other geometers; namely, in
\cite{L-X-S-J}, p. 126, we omitted the proof of Proposition 5.6.13 and
Proposition 5.6.15, stating that both proofs
 are very ``similar'' to the proofs of the foregoing
Propositions 5.6.12 and 5.6.14, resp.,  following exactly the
steps of these proofs.\footnote{X. Wang from the Australian National
University, another expert in the field, wrote a personal letter to
the third author; he complained  that our monograph \cite{L-X-S-J}
does not contain the proofs of the auxiliary Proposition 5.6.13 and
Proposition 5.6.15, resp.}  We decided to include some more details
here as {\it exemplary demonstration of the real affine techniques}.
For the convenience of the reader we give, where needed, precise
references to \cite{L-X-S-J} for a parallel reading, and guide the
reader to extended results.

Let us give a more precise description of the type of equation we
are going to investigate here: We study {\it Bernstein Properties}
of a nonlinear, fourth order partial differential equation for a
convex function $f$ on a convex domain $\Omega \subset
\mathbb{R}^{n};$ this equation can be written as
$$\sum_{i,j=1}^n F^{ij} w_{ij} = 0,
\;\;\;w :=\left[\det\left(\tfrac{\partial^{2}f}{\partial
x_{i}\partial x_{j}}\right)\right]^{a},\eqno(1.1)$$ where \enspace
$a\neq 0$ \enspace is a real constant and where $(F^{ij})$ \enspace
denotes \enspace the cofactor \enspace matrix of the Hessian
 matrix $(f_{ij}): =\left(\tfrac{\partial^{2}f}{\partial
x_{i}\partial x_{j}}\right) $. \vskip 0.1in \noindent Denote by
$\Delta$ and $\|\cdot\|$ the Laplacian and the tensor norm with
respect to the {\it Calabi metric} $\; G:= \sum f_{ij}dx_idx_j,\;$
respectively. Set
$$\rho:= \left[\det(f_{ij})\right]^{-\frac{1}{n+2}}
\;\; \;\; \text{and} \;\; \;\; \Phi : =\frac{1}{\rho^2} \cdot
\|\text{grad}\, \rho \|^2.
$$ In terms of the Calabi metric the PDE (1.1) can be rewritten as
$$\Delta \rho  = - \beta \; \tfrac{\|\text{grad}\  \rho \|^2}{\rho},
\eqno(1.2)$$ where $$\beta := -\,\tfrac{(n+2)(2a+1)+2}{2}.$$ The
\enspace PDE (1.1) \enspace appears in different geometric problems,
namely:
\begin{itemize}
\item {\it Chern's Conjecture:} \enspace In case when
$a = -\tfrac{n+1}{n+2}$, the PDE (1.1) is
the equation for affine maximal hypersurfaces; this is the
equation that is related to Chern's Conjecture cited above;
for Chern's Conjecture we have $n = 2$ and $\beta = 0.$
\item {\it Abreu equation:} \enspace  In  case that
$a = -1$ the PDE (1.1)  appears in the study of the differential
geometry of toric varieties (see \cite{A}, \cite{D}).
\end{itemize}

{\bf Remark.} \enspace Chern's completeness assumption is given by
the fact that the function $f$ is assumed to be defined for all $x
\in \mathbb{R}^n$; this completeness assumption is  $\;$  affinely
invariant, $\;$ in $\;$ affine hypersurface  theory $\;$ it is
called \\{\it Euclidean completeness.} \enspace Instead Calabi
assumed that the {\it Blaschke metric}  of the affine maximal
surface should be complete; this completeness assumption is
called {\it affine   completeness}.\\

About the Bernstein property for  the Abreu equation, A.-M. Li and
F. Jia posed the following conjecture in [L-J-1]: \vskip 0.1in
\noindent {\bf Conjecture \enspace of  Li-Jia.} \enspace  {\it Let
$f$ be a smooth, strictly convex function defined for all
$x\in\mathbb{R}^{n}.$ Assume that $f$ satisfies the Abreu equation
$$\sum_{i,j=1}^n F^{ij} w_{ij} = 0,
\;\;\;w: =\left[\det\left(\tfrac{\partial^{2}f}{\partial
x_{i}\partial x_{j}}\right)\right]^{-1}.$$ Then $f$ must be a
quadratic polynomial.}

\vskip 0.1in \noindent In [L-J-1] Li and Jia proved the following
theorem:

\vskip 0.1in \noindent {\bf Theorem 1.}\enspace (Theorem 5.6.2 in
\cite{L-X-S-J}).  \enspace  {\it Let $f$ be a smooth, strictly
convex function defined for all $(x_{1},x_{2})\in\mathbb{R}^2$. If
$f$ satisfies the PDE
$$\sum_{i,j=1}^2 F^{ij} w_{ij} = 0,
\;\;\;w :=\left[\det \left(\tfrac{\partial^{2}f}{\partial
x_{i}\partial x_{j}}\right)\right]^{a}$$ with $a\leq
-\frac{3}{4}$, then $f$ must be a quadratic polynomial.}
\vskip 0.1in \noindent {\bf Remark.}\\(1) \enspace When
$a = -\frac{3}{4}$,  Theorem 1  gives a new proof for
Chern's  conjecture on complete affine maximal surfaces.
That is a particular topic of this paper.\\
(2) \enspace When $a = - 1$,  Theorem 1 affirmatively solves the
above Conjecture of Li-Jia  for $n=2$.\\
(3) \enspace In [T-W-2], \enspace  N. Trudinger and X. Wang proved that the
global solution of the PDE (1.1) with $a>0$ on $\mathbb{R}^2$
must be a quadratic polynomial.\\
(4)  \enspace When $a=0$ and  $n=2$, the above PDE (1.2) reads
$\;\;\Delta \rho = 3\,\tfrac{\|\text{grad}\; \rho\|^2}{\rho}$; \\
this  is equivalent to
$$\sum_{i,j=1}^2 f^{ij}\tfrac{\partial^{2}}{\partial x_{i}\partial
x_{j}}[\ln\det(f_{kl})]=0.\eqno(1.3)$$ (5) \enspace  The global
solution of the PDE (1.3) on $\mathbb{R}^2$  is not unique. In fact,
the following are different global solutions of (1.3):
$$f(x_1,x_2)= \sum _{i=1}^2 x_i^2,\quad\quad  \quad \quad
h(x_1,x_2)= e^{x_1} + x_2^2.$$

In our foregoing sketch we indicated tools from analysis and
geometry by catchphrases  of experts.
In the following  short list we give some
hints and references:
\begin{itemize}
\item {\it Caffarelli-Schauder estimates:}\\
\cite{CAF} gives interior a priori estimates for solutions of
perturbations of   Monge-Amp\`{e}re equations; \cite{CAF-1}
treats a localization property of convex viscosity solutions of
certain Monge-Amp\`{e}re inequalities.
\item {\it Caffarelli-Guti\'{e}rrez theory:} \\
The paper \cite{CG} treats  the properties of the solutions of a
linearized Monge-Amp\`{e}re  equation.
\item {\it Affine blow-up analysis:} \\
One can apply a series of affine transformations to rescale the
domain and the function such that the geometric invariant to be
estimated has a better behavior. Such a  procedure is called {\it  affine
blow-up analysis}. In this process we often use a Lemma of Hofer.
Affine blow-up analysis is a very powerful tool to estimate affine
invariants in various circumstances, see \cite{L-J-2}, \cite{L-J-3},
\cite{L-J-4}, \cite{C-L-S-3}.
\item {\it Tools from the theory of convex bodies\,:}\\
A Theorem of Alexandrov-Pogorelov-Heinz gives sufficient conditions
for the local convexity of a surface; such conditions are of
importance for the solution of certain elliptic Monge-Amp\`{e}re
equations; see \cite{H} and \cite{BU},  p. 35.
\item {\it Affine and relative hypersurface theories:}\\
The extrinsic curvature functions in the affine theories are of
fourth order, in particular the Euler-Lagrange equation for the {\it
Affine Bernstein Conjecture}\,; the invariants that characterize
quadrics play an important role for the solution of certain
Monge-Amp\`{e}re equations; for details see the monographs
\cite{L-S-Z} and  \cite{S-S-V}.
\end{itemize}

In \cite{L-X-S-J} we presented an  analytic proof of  Theorem 1. In
this
note we outline a proof with emphasis on the proof of Chern's conjecture.\\

\section{ Proof of Theorem 1.}

Let $C>0$ be given a real constant and a convex domain $\Omega
\subset \mathbb{R}^n$; denote by $\mathcal {S}(\Omega,C)$ the class
of strictly convex $C^\infty$-functions $f$, defined on $\Omega$,
such that
$$\inf_{\Omega}f(x)= 0, \quad f(x) = C \quad \text{on}
\quad \partial \Omega.$$ We use the notion of a {\it normalized
convex set} in the sense stated on p. 48 in \cite{G}. To prove
Theorem 1  we need the following Lemma:

\vskip 0.1in \noindent {\bf Lemma 2.1.} \enspace (Lemma 5.6.16 in
\cite{L-X-S-J}). \enspace {\it Let $\Omega_k\subset\mathbb{R}^2$ be
a sequence of smooth normalized convex domains, converging to a
convex domain $\Omega $, and let $f^{(k)}\in
\mathcal{S}(\Omega_k,C)$ with $f^k(q^k)=0$. Assume that the
functions $f^{(k)}$ satisfy the PDE (1.2) with $\beta\geq0$. Then
there exists a subsequence $f^{(i_{\ell})}$ that locally uniformly
converges to a convex function $f\in C^{0}(\Omega)$ with distance
$dist(p_o,\partial \Omega )>0$, where $p_o$ is the point such that
$f(p_o)=0$. Moreover, there is an open neighborhood $N$ of $p_o$
such that $f^{(i_{\ell})}$ converges to $f$, and also all their
derivatives converge, therefore $f$ is smooth and strictly convex in
$N$.}

\vskip 0.1in \noindent We will sketch the proof of this Lemma in
Section 3.

\vskip 0.1in \noindent {\bf Proof of Theorem 1.}  \enspace Let
$x:M\rightarrow \mathbb{R}^{3}$ be a locally strongly convex
surface, given \enspace as \enspace graph of a smooth, strictly
\enspace convex \enspace function \enspace $f,$ \enspace defined for
all $(x_{1},x_{2})\in \mathbb{R}^{2}$. Assume that $f$ satisfies the
PDE (1.2) with $\beta\geq 0$. Given any $p\in M$, - adding a linear
function if necessary, - we may assume that
$$f(p) = 0,\;\;\;\;\;\tfrac{\partial f}{\partial x_{i}}
(p)=0,\;\;\; i=1,2.$$ Choose a sequence $\{C_{k}\}$ of real positive numbers
such that $C_{k}\rightarrow\infty$ as $k\rightarrow\infty$. Then,
for any $C_{k},$ the \emph{section}
$$S_{f}(p,C_k) := \left\{(x_{1},x_{2})\in \mathbb{R}^{2}\;|\;\;
 f(x_1,x_{2}) < C_{k}\right\}$$
is a bounded convex domain in $\mathbb{R}^{2}$. It is well-known
(see section 1.8 in \cite{G}) that there exists a unique ellipsoid
$E_{k}$ which attains the minimum volume among all  ellipsoids that
contain $S_{f}(p,C_k)$ and that are centered at the center of mass
of $S_{f}(p,C_k)$ such that
$$2^{-\frac{3}{2}} E_{k} \subset S_{f}(p,C_k) \subset E_{k}\,.$$
Let $T_{k}$ be an affine transformation such that $$
T_{k}(E_{k})=B_1(0)=\left\{(x_{1},x_{2})\in\mathbb{R}^{2}
\;\;|\;\;(x_{1})^{2}+ (x_{2})^{2}< 1\right\}.$$ Define the functions
$$f^{(k)}(x):=\tfrac{1}{C_{k}}\;f(T^{-1}_{k}x)    .$$
Then
$$B_{2^{-\frac{3}{2}}}(0)\subset\Omega_{k}\subset B_1(0),$$ where
$$\Omega_{k}=\left\{(x_{1},x_{2})\in \mathbb{R}^{2}\;\;|\;\;
 f^{(k)}(x_{1},x_{2}) < 1 \right\}.$$
Taking subsequences, we may assume that $\{\Omega_k\}$ converges to
a convex domain $\Omega$ and $\{f^{(k)}\}$ converges to a convex
function $f^{\infty}$, locally uniformly in $\Omega$. By Lemma 2.1
the function $f^{\infty}$ is smooth and strictly convex in a
neighborhood of $T_k(p)\in\Omega$, for $k$ large enough. It follows
that the functions
$\Phi^{(k)}(T_k(p))$ are uniformly bounded.\\
Assume that $\Phi (p)\ne 0$;  a direct calculation gives
$$\Phi^{(k)}(T_k(p))= C_{k}\Phi (p) \rightarrow \infty,$$
thus we get a contradiction, and thus $\Phi (p) = 0$. Since $p$ is
arbitrary, we have $\Phi = 0$ everywhere on $M$. It follows that
$\det (f_{ij}) = const.$  The Theorem of J\"orgens (\cite{L-X-S-J},
p. 59) implies that $f$ must be a quadratic polynomial. This
completes the proof of Theorem 1. \hfill $\blacksquare$

\section{ Proof of Lemma 2.1.}

In this section we sketch two different proofs  of Lemma 2.1. The
first  one uses the theory  of convex bodies  combined with  the
theory of Caffarelli-Guti\'{e}rrez, while  the second one  uses pure
analysis. Both proofs need the following estimates for determinants.
Here and later we separate the cases $\beta > 0$ \enspace and
\enspace $\beta = 0.$
\subsection{ Estimates for  determinants}

\vskip 0.1in \noindent  We restrict to dimension $n = 2.$  For $\beta>0$ we
have:

\vskip 0.1in \noindent {\bf Lemma 3.1.} \enspace (Lemma 5.6.9 in
\cite{L-X-S-J}). \enspace {\it Let $0 < C \in \mathbb{R}$ and
$\Omega$ be a normalized convex domain with  $0\in \Omega$ as center
of $\Omega$.  Assume that $f\in\mathcal {S}(\Omega,C)$ and $f$
satisfies the PDE (1.2) with $\beta>0$. Let $u$ be the Legendre
transformation function of $f$ relative to $0$. Then the following
estimate holds:
$$\tfrac{1}{(d+u)^{4}}\cdot \det(f_{ij}) \leq b_0
\qquad\mbox\ for\;\; \;x\in S_f(C')= \left\{x\in\Omega\;|\;f(x)\leq
C'\right\},$$ where $0 < C'< C$, and $b_0$ is a constant depending
only on $\frac{C'}{C}$, $C$, and $\beta $.} \vskip 0.1in \noindent
The affine maximal surface equation is affinely invariant; thus, for
$\beta = 0$, we can choose a  coordinate system
$\;\;x_1,x_2,x_3\;\;$ such that the graph function \\
$x_3=f(x_1,x_2)$  satisfies  $\|grad\, f \|_E \leq C_0$ for some
constant $C_0$ in a section; here and later $\|\cdot\|_E$ denotes
the norm with respect to the canonical Euclidean metric of
$\mathbb{R}^2$. From Lemma 5.6.4 in \cite{L-X-S-J} we have

\vskip 0.1in \noindent {\bf Lemma 3.2.} \enspace  {\it Let $C$ and
$\Omega$  be defined as before. Assume that $f\in\mathcal
{S}(\Omega,C)$ and $f$ satisfies the PDE (1.2) with $\beta=0$. If
$$\|grad\, f\|_E \leq C_0$$ for some constant $C_0$, then the following estimate holds:
$$\det(f_{ij})\leq b_0  \qquad\mbox\
for\;\;\; x\in S_f(C'),$$ where $b_0$ is a constant depending only
on $0 < C$, $\;0 <  \frac{C'}{C} < 1\;$ and $\;C_0$.}

\vskip 0.1in \noindent {\bf Lemma 3.3.} \enspace  (Lemma 5.6.6 in
\cite{L-X-S-J}). \enspace {\it Let $f$ be a smooth, strictly convex
function defined on a bounded convex domain $\Omega \subset
\mathbb{R}^{2}$, satisfying the PDE (1.2). Denote by $\Omega^*$ the
Legendre transformation domain of $f$. Let $\Omega^{\prime *}$ be an
arbitrary subdomain of $\Omega^*$ with dist$\,(\Omega^{\prime
*},\partial\Omega^*)>0$. Then the following estimate holds:
$$\det(f_{ij})\geq b_0\quad \mbox\ for\quad
\xi\in\Omega^{\prime *},$$ where $b_0$ is a constant depending only
on dist$(\,\Omega^{\prime
*},\partial\Omega^*)$, $diam(\Omega )$, $diam(\Omega^*)$ and
$\beta \in \mathbb{R}$.}
\vskip 0.1in \noindent

One can use the Lemmata  3.1 - 3.3, the theory of convex bodies, and
finally the Caffarelli-Guti\'{e}rrez theory to prove Lemma 2.1. This
was carried out in [L-J-3] for $\beta =0$ (see also section 6.1 in
\cite{L-X-S-J}). In a remark in [L-J-1],  Li and Jia have pointed
out that this method works for $\beta>0$ (see also p. 129 in
\cite{L-X-S-J}); but they did not  give details. Here we give a
proof.
\vskip 0.1in \noindent {\bf First Proof of Lemma 2.1.} \vskip 0.1in
\noindent This proof uses the Caffarelli-Guti\'{e}rrez theory and
tools from the theory of convex bodies. \\
{\bf Case
$\beta > 0$.} \enspace First we consider the case
$\beta > 0$. From Lemma 3.1, for the function $f^{(k)}$ in Lemma 2.1,
we have the following estimate  in the section
$S_{f^{(k)}}(q^{k},\tfrac{C}{2}):$
$$\frac{\det(f^{(k)}_{ij})}{(d+u^{(k)})^4}\leq d_0.\eqno(3.1)$$
Choose the radius $R$ for  a Euclidean ball  such that
 $$ B_{R}(q^{k})\supset\Omega^{(k)}\supset
S_{f^{(k)}}(q^{k},\tfrac{C}{2}); \quad \quad \text{then \, the \,
relation }\quad \quad B^\ast_r(0)\subset (\Omega^{(k)})^\ast
$$
is satisfied for the Legendre transformation domains, where
$r=\frac{C}{4 R}$. Moreover, restricting to  $B^\ast_r(0)$, we have
the inequality
$$-\frac{C}{2R}-C\leq u^{(k)}\leq \frac{C}{2R}.$$
Therefore,
\begin{enumerate}
\item $u^{(k)}$ locally converges to a convex function $u^\infty$
on $B^\ast_r(0)$;
\item from standard relations of the Legendre transformation
(see  section 1.4 in \cite{L-X-S-J}) and from (3.1)  we know that
$\det(u_{ij}^{(k)})=[\det (f_{ij}^{(k)})]^{-1}>d_1$ is bounded from
below on $B^\ast_r(0)$.
\end{enumerate}
Then the Alexandrov-Pogorelov-Heinz Theorem implies that $u^\infty$
is strictly convex in a small neighborhood of $(0,0)\in
B^\ast_r(0)$.
Adding a linear function if necessary, we may assume that
$$u^\infty(0) = 0,\;\;\;\;\;\tfrac{\partial u^\infty}{\partial \xi_{i}}
(0)=0,\;\;\;i=1,2.$$
Thus  there exists  a real constant $h_0>0$ such that the section
$$S_{u^\infty}(h_0):=\{\xi\in B^\ast_{\frac{r}{2}}(0)\;\;|\;\;
u^\infty\leq h_0\}$$ is compact in $B^\ast_{\frac{r}{2}}(0)$. Then, in
this section $S_{u^\infty}(h_0)$, we have a lower bound for
$\det(u_{ij}^{(k)})$ and, by Lemma 3.3, we also have an upper bound.
 Now  we use the Caffarelli-Guti\'{e}rrez theory and
the Caffarelli-Schauder estimate to conclude that $\{u^{(k)}\}$
smoothly converges to $u^{\infty}$ in $\{\xi\; |\; u(\xi)\leq
h_0\}$. Therefore, $u^{\infty}$ is  a smooth and strictly convex
function in an open neighborhood of (0,\;0) in
$B^*_{\frac{r}{2}}(0)$. \vskip 0.1in \noindent {\bf Case $\beta =
0$}. \enspace Next we consider the case $\beta = 0$ ( see
\cite{L-X-S-J}, p.129 and   p. 153-159). Denote
$D:=\{x\;|\;f(x)=0\}$. Again we discuss  two subcases.

\vskip 0.1in \noindent {\bf(i)} If $D\bigcap\partial \Omega =
\emptyset$ then there is a constant $h>0$ such that the level set
satisfies $\bar{S}_f(p_o,h)\subset \Omega$, and thus  we have a
uniform estimate for $\|grad \,f^{(k)}\|_E$ in
$\bar{S}_{f^{(k)}}(q^k,h)$. From Lemma 3.2  it follows that there is
a uniform estimate for $\frac{1}{\rho}$ in
$\bar{S}_{f^{(k)}}(q^k,\frac{h}{2}).$ Then we use the same argument
as in the case $\beta>0$ to conclude that $\{u^{(k)}\}$ smoothly
converges to $u^{\infty}$ in $\{\xi\; |\; u(\xi)\leq h_0\}$.
Therefore  $u^{\infty}$ is  a smooth and strictly convex function in
an open neighborhood of (0,\;0) in $B^*_{\frac{r}{2}}(0)$.

\vskip 0.1in \noindent{\bf(ii)} If  $D\bigcap \partial \Omega \ne
\emptyset$, let $p\in D\bigcap\partial \Omega $. Since the PDE (1.2)
with $\beta = 0$ is equiaffinely  invariant, we may choose a new
coordinate system such that the term $\|grad \,f^{(k)}\|_E$ is
uniformly bounded in $\bar{S}_{f^{(k)}}(p,h)$. Then the same
argument shows that $f$ is smooth in a neighborhood of $p$, and we
get a contradiction. This excludes the case $D\bigcap
\partial \Omega \ne \emptyset$.\\
Lemma 2.1 is proved. \hfill
$\blacksquare$
\vskip 0.1in \noindent
{\bf 3.2 Estimates for the third derivatives
and for $\sum f_{ii}$.}
\vskip 0.1in  It is  interesting to give a purely analytic proof of
Lemma 2.1, without using tools from the theory of convex bodies.
Here we apply the {\it real affine technique} that we mentioned in
the Introduction.
The key point of the proof of Lemma 2.1 are the estimates for the
third order derivatives and for $\sum f_{ii}$. We introduce the
following notations:
$$\mathcal {A} : = \max_{ \Omega}\left\{\exp \left \{-\tfrac{m}{C-f} \right \}
\tfrac{\Phi}{\rho^{\alpha}(d+u)^{\alpha}}\right\},$$
$$\mathcal {D}:= \max_{\Omega}\left\{\exp \left \{-\tfrac{m}
{C-f}+ K \right \}\tfrac{g^2\,
\|\,grad\;f\|^2}{\rho^{\alpha}(d+u)^{\alpha}}\right\},$$ where
$$K := \tfrac{N}{\mathcal {A}} \exp \left \{-\tfrac{m}{C-f}\right \}\tfrac{\Phi}
{\rho^{\alpha}(d+u)^{\alpha}},$$ and $m$, $\alpha $ and $N$ are
positive real constants. As before we separate the cases $\beta > 0$
\enspace and \enspace $\beta = 0.$

\vskip 0.1in \noindent {\bf Lemma 3.4.} \enspace (Lemma 5.6.10 in
\cite{L-X-S-J}). \enspace {\it Assume that
$\Omega\subset\mathbb{R}^2$ is a normalized domain and $f \in
\mathcal {S}(\Omega,C)$ satisfies the PDE (1.2) with $\beta>0$, and
that there exists a constant $b>0$ such that in
$\Omega$:$$\tfrac{1}{\rho(d+u)}<b.$$ Then there are constants
$\alpha>0$, $N$ and $m$ such that
$$ \mathcal {A}\leq \max\{d_1, \tfrac{4}{\alpha \beta}\mathcal {D}\},$$
where $d_1$ is a real constant depending only on $\alpha$, $C$, $b$ and
$\beta$.}

\vskip 0.1in \noindent {\bf Lemma 3.5.} \enspace  (Lemma 5.6.11 in
\cite{L-X-S-J}). \enspace {\it  Let $\Omega\subset\mathbb{R}^2$ be a
normalized convex domain and $f \in \mathcal {S}(\Omega,C)$,
satisfying the PDE (1.2) with $\beta>0$. Assume that there exists a
constant $b>0$ such that in $\Omega$:
$$\tfrac{1}{\rho(d+u)}<b.$$
Then there exist real constants $\alpha >0$, $N$ and $m$ such that
$$\mathcal {A}\leq d_2,\;\;\;\;\;\;\mathcal {D}\leq d_2$$
for some real constant $d_2>0$ that depends only on $\alpha$, $b$,
$\beta$ and $C$.}

\vskip 0.1in As a corollary of  Lemma 3.4 and Lemma 3.5  we get:

\vskip 0.1in \noindent{\bf Proposition 3.6.} \enspace (Proposition
5.6.12 in \cite{L-X-S-J}). \enspace   {\it Let
$\Omega\subset\mathbb{R}^2$ be a normalized convex domain and
$0\in\Omega $ be the  center of $\Omega$. Let $f$ be a strictly
convex $C^{\infty}$ function defined on $\Omega$. Assume that
$$\inf_\Omega f=0,\quad\quad f=C>0\;\;\;on\;\;\partial \Omega,$$
and that $f$ satisfies the PDE (1.2)  with $\beta
>0$. Then there exists a constant $\alpha >0$ such that,
on $\Omega_{\frac{C}{2}}:=\{x\in\Omega\;|\; f(x)<\frac{C}{2}\}$,
there is a joint upper bound
$$\tfrac{\Phi}{\rho^\alpha(d+u)^\alpha}\leq d_3,\;\;\;
\tfrac{\|\,grad\; f\|^2 }{\rho^\alpha(d+u)^\alpha} \leq d_3$$ for
some constant $d_3>0$ that depends only on $\beta$ and $C$.}

\vskip 0.1in \noindent {\bf Proposition 3.7.}\enspace (Proposition
5.6.14 in \cite{L-X-S-J}).\enspace {\it Let
$\Omega\subset\mathbb{R}^2$ be a normalized convex domain. Let $f
\in\mathcal {S}(\Omega,C)$ be a smooth and strictly convex function
defined in $\Omega$, which satisfies the PDE (1.2) with $\beta >0$.
Assume that there are constants $d_3>0$ and $\alpha>1$ such that, in
$\Omega$,
$$\tfrac{\Phi}{\rho^\alpha(d+u)^\alpha}<d_3,\quad \quad \quad
\tfrac{1}{\rho(d+u)}<d_3.$$ Then there exists a constant $d_5>0$,
depending only on $\alpha$, $\beta$, $d_3$ and $C$, such that, on $\Omega$:
$$
\exp \left\{-\tfrac{32(2+d_3)C}{C-f }\right\}\tfrac{\sum f_{ii}
}{\rho^{\alpha }(d + u)^{\alpha+2}}\leq d_5 \,.
$$}

For affine maximal surfaces we have $\beta = 0$ in the PDE (1.2). We
already stated  above  that,  for the equation of affine maximal
surfaces, we may choose an appropriate  coordinate system  reducing
our problem to the case that $\|\text{grad}\, f\|_E$ \enspace is
bounded above, and  then proving that \enspace $\frac{1}{\rho}\leq
b$ \enspace for a certain constant $b$. In this case we have \vskip
0.1in \noindent{\bf Proposition 3.8.} (Proposition 5.6.13 in
\cite{L-X-S-J}). \enspace  {\it Let $\Omega\subset\mathbb{R}^2$ be a
normalized convex domain and $0\in\Omega $ be the  center of
$\Omega$. Let $f$ be a strictly convex $C^{\infty}$ function defined
on $\Omega$. Assume that
$$\inf_\Omega f=0,\quad\quad f=C>0\;\;\;on\;\;\partial \Omega,$$
and that $f$ satisfies  the PDE (1.2) with $\beta=0$; moreover
assume  that there is a constant $b>0$ such that, in $\Omega$:
$$\tfrac{1}{\rho} \leq b.$$ Then there exists a real constant
$\alpha >0$ such that the following estimates hold on
$\Omega_{\frac{C}{2}}$:
$$\tfrac{\Phi}{\rho^\alpha}\leq d_4,\;\;\;\;\;\;
\tfrac{\|\,\text{grad}\; f\|^2 }{\rho^\alpha} \leq d_4$$ for some
real constant $d_4>0$ that depends  only on $\alpha$, $b$ and $C$.}

\vskip 0.1in \noindent {\bf Proposition 3.9.} \enspace (Proposition
5.6.15 in \cite{L-X-S-J}). \enspace {\it Let $x_{3}=f(x_1,x_{2})$ be
a smooth and strictly convex function defined on a normalized convex
domain $\Omega \subset \mathbb{R}^2$, \enspace which satisfies the
equation \enspace (1.2) \enspace with $\beta = 0$. \enspace Assume
that there exist constants $\alpha \geq 0$ and $d_4\geq 0$ and a
function $f\in \mathcal {S}(\Omega,C)$  such that:
$$\tfrac{\Phi }{\rho^{\alpha}}\leq d_4,\;\;\;\; \tfrac{1}{\rho}\leq
d_4$$ on $\bar{\Omega}$. Then there is a real constant $d_5
>0$, depending only on $\alpha $, $d_4$ and $C$, such that on $\Omega$
$$\exp \left\{-\tfrac{32(2+d_4)C}{C-f }\right\}\tfrac{\sum f_{ii}
}{\rho^{\alpha }(d+u)^{\alpha+2}}\leq d_5\, .$$ }
{\bf Remark.} It might be helpful for the understanding of the
foregoing proofs to realize certain differences
in the proofs, namely:\\
{\bf (1)}  For $\beta > 0$ we have the estimate of the determinant
$$\tfrac{1}{(d+u)^{4}}\cdot \det(f_{ij}) \leq b_0
\qquad\mbox\ for\;\; \;x\in S_f(C')= \left\{x\in\Omega\;|\;f(x)\leq
C'\right\}$$ from Lemma 3.1, while for $\beta=0$ we don't have  such
an estimate. But in the latter  case the PDE (1.2) is affinely
invariant, thus we can choose affine coordinates $x_1, x_2,x_3$ such
that $x_3=f(x_1,x_2)$ satisfies $\|grad\; f\|_E\leq C_0$ for some
constant $C_0>0$; as a consequence,  from Lemma 3.2,  we then have
the estimate
$$\det(f_{ij}) \leq b_0  \quad \text{for} \;\; x\in S_f(C').$$
{\bf (2)}  The assumptions in  Proposition 3.6 and Proposition 3.8
are different, therefore we can not prove Proposition 3.8 as
Corollary from Proposition 3.6,  just considering the limit
$\beta\rightarrow 0$.

\vskip 0.1in \noindent In \cite{L-X-S-J} and [L-J-1] we proved in
detail both Lemmas 5.6.10 and 5.6.11, thus  also Proposition 5.6.12
in \cite{L-X-S-J} has been proved in all details. As the reader can
realize, the proofs of Proposition 5.6.13 and 5.6.15 are, step by
step,   similar to the proofs of Propositions 5.6.12 and 5.6.14,
including the method, the  estimates, and even the details of
computation; therefore  we omitted the  proofs of Propositions
5.6.13 and 5.6.15 in both publications, [L-J-1] and \cite{L-X-S-J}.
\vskip 0.1in \noindent {\bf Second Proof of Lemma 2.1.}

\vskip 0.1in \noindent This proof is taken from \cite{L-X-S-J}, p.
129; it uses pure analysis.\\
{\bf Proof.} \enspace We treat the two cases \enspace $\beta >0$
\enspace and \enspace $\beta = 0$ \enspace separately.

\vskip 0.1in \noindent {\bf Case $\beta >0$}. \enspace Let
$0\in\Omega_{k} $ be the center of $\Omega_{k}$ and $u^ {(k)}$ the
Legendre transformation function of $f^ {(k)}$ relative to 0.\\
To simplify the notations we will use $f^{(k)}$ to denote
$f^{(i_{\ell})}$. Lemma 3.1 and Propositions 3.6 and  3.7 imply the
following uniform estimates:
$$\tfrac{\Phi^{(k)} }{\rho^{\alpha}(d + u^{(k)})^{\alpha}}
\leq d_6,\;\;\;\; \tfrac{1}{\rho^{\alpha}(d + u^{(k)})^{\alpha}}\leq
d_6, \;\;\; \tfrac{\sum f^{(k)}_{ii} }{\rho^{\alpha }(d +
u^{(k)})^{\alpha+2}}\leq d_6 $$ in
$$S_{f^{(k)}}(q^k,\tfrac{C}{2}):=\{x\in\Omega_k\,|\,
 f^{(k)}<\tfrac{C}{2}\},$$ where $d_6$ is a positive constant
depending only on $\beta$ and $C$. We may assume that $q_k$
converges to $p_o$. Let $B_R(q^k)$ be a Euclidean ball with
$\Omega \subset B_{\frac{R}{2}}(q^k)$. Then the Legendre
transformation domain of $\Omega$ satisfies
$B^*_{\delta}(0)\subset\Omega^*$, where $\delta = \frac{C}{2R}$ and
$B^*_{\delta}(0) = \{\, \xi\;\, |\;\; \xi_1^2 + \xi_2^2 < \delta^2\}.$
Lemma 3.3 gives
$$\det(f_{ij})\geq b_3$$ for $\xi\in B^*_{\frac{\delta}{2}}(0),$
where $b_3$ is a constant depending only on $C$ and $\beta$.
Restricting to $B^*_{\delta}(0)$, we have
$$-\tfrac{C}{2R}-C\leq u^{(k)} =
\sum \xi_i x_i - f^{(k)} \leq \tfrac{C}{2R}\,.$$ Therefore  the
sequence  $u^{(k)}$ locally uniformly converges to a convex function
$u^{\infty}$ in $B^*_{\delta}(0)$, and there are constants \enspace
$0<\lambda \leq \Lambda < \infty$  \enspace such that the following
estimates hold in $B^*_{\frac{\delta}{2}}(0)$
$$\lambda \leq\lambda_i^{(k)}
\leq \Lambda,\qquad \hbox{for}\;\;\;i=1,2,\cdots,n, \;\;\;\;\,
k=1,2,\cdots$$ where $\lambda_1^{(k)},\cdots,\lambda_n^{(k)}$ denote
the eigenvalues of the matrix $(f^{(k)}_{ij})$. Standard elliptic
estimates finally prove the assertion of Lemma 2.1  in the case
$\beta>0$.

\vskip 0.1in \noindent {\bf Case $\beta = 0$}.\enspace Denote
$D:=\{x\;|\;f(x)=0\}$. \enspace Again we have two subcases.  \vskip
0.1in \noindent{\bf(i)} If $D\bigcap\partial \Omega = \emptyset$
then there is a constant $h>0$ such that the level set satisfies
$\bar{S}_f(p_o,h)\subset \Omega$, thus  we have a uniform estimate
for $\|grad \,f^{(k)}\|_E$ in $\bar{S}_{f^{(k)}}(q^k,h)$. From Lemma
3.2  it follows that there is a uniform estimate for
$\frac{1}{\rho}$ in $\bar{S}_f^{(k)}(q^k,\frac{h}{2}).$ Now  we use
Propositions 3.8 and 3.9 and the same argument as in the case $\beta
> 0$ in the second proof to complete the proof.

\vskip 0.1in \noindent {\bf(ii)} In case $D\bigcap \partial \Omega
\ne \emptyset$, let $p\in D\bigcap\partial \Omega $. Since the PDE
(1.2) with $\beta = 0$ is equiaffinely  invariant, we may choose a
new coordinate system such that $\|grad \,f^{(k)}\|_E$ is uniformly
bounded in $\bar{S}_{f^{(k)}}(p,h)$. Then the same argument as in
(i) above shows that $f$ is smooth in a neighborhood of $p$, and we
get a contradiction. This excludes the case $D\bigcap
\partial \Omega \ne \emptyset$ and thus completes the proof of
Lemma 2.1. \hfill $\blacksquare$
\section{ Proof of Proposition 3.8}
When formulating Proposition 3.8, we already have stated
that this Proposition is exactly Proposition 5.6.13 in
\cite{L-X-S-J}, but its proof was omitted in our  monograph for
the reasons given above.\\
To emphasize the similarity of both proofs we exceptionally use a
notation that we introduced in  section 3.2 above for similar, but
not identical  expressions; this might be allowed as the proof of
both Lemmas, namely 3.4 and 3.5, has been finished. Now we define:
\begin{alignat*}{1} \mathcal {A} : &= \max_{ \Omega}\left
\{\exp \left \{-\tfrac{m}{C-f} \right \}
\tfrac{\Phi}{\rho^{\alpha} }\right\},\\
\mathcal {D}:&= \max_{\Omega}\left\{\exp \left \{-\tfrac{m} {C-f}+ K
\right \}\tfrac{g^2\, \|\,\text{grad}\; f\|^2}{\rho^{\alpha}
}\right\},\end{alignat*} where
$$K := \tfrac{N}{\mathcal {A}}\, \exp
\left \{-\tfrac{m}{C-f}\right \}\tfrac{\Phi} {\rho^{\alpha} },$$ and
$m$, $\alpha $ and $N$ are positive constants to be determined
later. \\
It is important for the understanding of the following proof to
realize that we use two different differential inequalities for the
expression $\tfrac{\Delta \Phi}{\Phi}$, namely one coming from the
affine maximal surface equation, the other from the study of the
function $L$ below. The steps of the proof of Proposition 3.8 are
the following: \vskip 0.1in \noindent
\textbf{Step 1.} First we will prove that \enspace $\mathcal
{A}\leq\tfrac{4}{\alpha}\mathcal {D}.$  \enspace  The method is the same as in
the proof of Lemma 5.6.10 in  \cite{L-X-S-J}. Consider the function
$$L : = \exp \left \{-\tfrac{m}{C-f}\right \}\, \tfrac{\Phi}{\rho^{\alpha}}$$
defined on $\Omega$. From  the definition  of $L$ and its behavior
at the boundary, the function $L$ attains its supremum at some
interior point $p^*$. Differentiation gives  at $p^*$:
$$\tfrac{\Phi_{,i}}{\Phi} -  gf_{,i} - \alpha \tfrac{\rho_{,i}}{\rho} = 0,\eqno(4.1) $$
$$\tfrac{\Delta \Phi}{\Phi}-
\tfrac{\|\,\text{grad}\; \Phi \|^2}{\Phi^2}+ \alpha \Phi - g'\|\,
\text{grad} \,f \|^2 - g\Delta f\leq 0; \eqno(4.2)$$
 here in (4.2) and  later we use the abbreviation
$$g:=\tfrac{m}{(C-f)^2},\quad\quad g':=\tfrac{2m}{(C-f)^3}.$$
We choose $m=2C$. Then $g'\leq g^2$. Then the inequality in (4.2)
gives an upper bound for the expression
$\tfrac{\Delta \Phi}{\Phi}$.\\
On the other hand, the affine maximal surface equation (see
Proposition 4.5.2 in \cite{L-X-S-J}, in this case $n=2$ and
$\beta=\delta=0$) gives the following
 lower bound for $\tfrac{\Delta \Phi}{\Phi}$:
$$\Delta \Phi\geq
\tfrac{\|\text{grad}\; \Phi\|^2}{\Phi}+ \Phi^2.\eqno(4.3)$$
A combination of the inequalities (4.2), (4.3) and
another application of the   inequality of Schwarz  give
$$\exp \left \{-\tfrac{m}{C-f}\right\}\tfrac{\Phi}{\rho^{\alpha}} \leq
\tfrac{ 2 }{ \alpha }\exp \left \{-\tfrac{m}{C-f}\right
\}g^{2}\,\tfrac{\|\,\text{grad} \, f\|^2}{\rho^{\alpha}} +
\tfrac{4b^\alpha}{\alpha C}. $$   We may assume that, at $p^*$, $
\exp \left \{-\tfrac{m}{C-f}\right \}g^{2}\tfrac{\|\, \text{grad} \,
f\|^2}{\rho^{\alpha}} > \tfrac{2b^\alpha}{C}.$ Note that \\ $K(p^*)=
N.$ Then  Step 1 is proved.

\vskip 0.1in \noindent \textbf{Step 2.} We will prove that $\mathcal
{D} $ has an  upper bound.  Consider the function (see (5.6.6) in
\cite{L-X-S-J})
$$F: = \exp \left \{\tfrac{-m}{C-f} + \tau \right \}\,Q\,\|\, \text{grad}\;
h\,\|^2$$ defined on $\Omega$, and put
$$\tau:= K ,\;\;\; Q:=
\tfrac{g^2}{\rho^{\alpha} },\;\;\;h:= f.$$  Assume  that $F$ attains
its supremum at the point $q^*$. Choose a local orthonormal frame
field with respect to the Calabi metric near $q^*$ such that
$f_{,1}(q^*)= \|\,\text{grad}\,f\|$, and choose $\delta =
\tfrac{1}{10}$, \enspace $N>>10$. Then   a  calculation as in the
proof of Lemma 5.6.11 in \cite{L-X-S-J}    gives, at $q^*$:
$$\left(-gf_{,i}  + \tfrac{4}{C-f}f_{,i}
- \alpha \tfrac{\rho_{,i}}{\rho}+ K_{,i} \right) (f_{,1})^2 + 2\sum
f_{,j}f_{,ji} = 0,\eqno(4.4)$$
$$2 (f_{,11})^2 +  2 (f_{,12})^2 +
4\tfrac{\rho_{,11}}{\rho}(f_{,1})^2 + (\alpha - 388)\Phi (f_{,1})^2
- 328 + \Delta K(f_{,1})^2$$$$  - \left(g'(f_{,1})^2 + 2\left(g -
\tfrac{4}{C-f}\right)+ 2\left(g -
\tfrac{4}{C-f}\right)\tfrac{\rho_{,1}}{\rho}f_{,1}\right)
(f_{,1})^2\leq 0.\eqno(4.5)$$
In the following we calculate estimates
for three terms appearing in (4.5), namely the terms:\\
\vspace{2mm}
-  \enspace $(f_{,11})^2 + (f_{,12})^2$.\\
\vspace{2mm}
-  \enspace $\Delta K$.\\
\vspace{2mm}
-  \enspace $4\tfrac{\rho_{,11}}{\rho}(f_{,1})^2$, \, resp.\\
Again we emphasize that  these  steps and details of
our calculation
are similar to the steps and details
in case that $\beta > 0;$ see  pp. 123-126
in \cite{L-X-S-J}.
\vskip 0.1in \noindent {\bf Estimate for the term  \enspace $(f_{,11})^2 + (f_{,12})^2$.}
$$
2(f_{,11})^2 = \tfrac{1}{2}\left[gf_{,1} - \tfrac{4}{C-f}f_{,1} +
\alpha  \tfrac{\rho_{,1}}{\rho} - K_{,1} \right]^2(f_{,1})^2 $$$$
\geq \tfrac{3}{4N} \left[\left(g - \tfrac{4}{C-f}\right)f_{,1} +
\alpha  \tfrac{\rho_{,1}}{\rho}  \right]^2(f_{,1})^2 -
\tfrac{9}{10}\tfrac{(K_{,1})^2}{K}(f_{,1})^2,\eqno(4.6)$$
$$2(f_{,12})^2 \geq \tfrac{3}{4N}\alpha^2
(\tfrac{\rho_{,2}}{\rho} )^2(f_{,1})^2  -
\tfrac{9}{10}\tfrac{(K_{,2})^2}{K}(f_{,1})^2.\eqno(4.7)$$

\vskip 0.1in \noindent {\bf Estimate for the term \enspace  $\Delta K$.}
$$K_{,i}=K\left(\tfrac{\Phi_{,i}}{\Phi}-
\alpha\tfrac{\rho_{,i}}{\rho}-gf_{,i}\right),\eqno(4.8)$$
$$\Delta K \geq \tfrac{\|\text{grad}\; K\|^2}{K}-2Kg\;\tfrac{
\rho_{,1}f_{,1}}{\rho} - N \left(g'(f_{,1})^2+2g
\right).\eqno(4.9)$$

\vskip 0.1in \noindent {\bf Estimate for the term  \enspace \enspace
$4\,\tfrac{\rho_{,11}}{\rho}\,(f_{,1})^2$. }
\vskip 0.1in \noindent Note that $\beta=0$ now, by the same
estimates as in Lemma 5.6.11 ( \cite{L-X-S-J},  p. 124), we have
$$\sum\tfrac{(\rho_{,ij})^2}{\rho^{2}} \leq \tfrac{2((\rho_{,11})^2 +
(\rho_{,12})^2)}{\rho^2}\leq \tfrac{\sum (\Phi_{,i})^2}{\Phi} +
4\Phi^2.$$ It follows that $$ 4\tfrac{|\rho_{,11}|}{\rho} (f_{,1})^2
\leq  4\;
\tfrac{\left\|\text{grad}\;\Phi\right\|}{\sqrt{\Phi}}(f_{,1})^2+
8\;\Phi (f_{,1})^2\leq  8\Phi
(f_{,1})^2 + $$$$ +  4\sqrt{2}\sqrt{\Phi}\left(
\left[\sum\left(\tfrac{\Phi_{,i}}{\Phi}-gf_{,i}- \alpha
\tfrac{\rho_{,i}}{\rho}\right)^2\right]^{\frac{1}{2}}
 + \left[\sum\left(gf_{,i}+ \alpha\tfrac{\rho_{,i}}{\rho}
\right)^2\right]^{\frac{1}{2}}\right)(f_{,1})^2 . \eqno(4.10)$$ We
apply the inequality of Schwarz and (4.8) to get
$$
4\sqrt{2}\sqrt{\Phi}\left[\sum\left(\tfrac{\Phi_{,i}}{\Phi}
-gf_{,i} - \alpha
\tfrac{\rho_{,i}}{\rho}\right)^2\right]^{\frac{1}{2}} (f_{,1})^2 \leq  $$
$$
\leq \tfrac{1}{20}\tfrac{\|\text{grad}\; K\|^{2}}
{K}(f_{,1})^2+\tfrac{240}{N}\mathcal {A}\,\exp\left\{
\tfrac{m}{C-f}\right\}\rho^{\alpha}(f_{,1})^2,$$
$$4\sqrt{2}\sqrt{\Phi}\left[\sum\left(gf_{,i}+\alpha
\tfrac{\rho_{,i}}{\rho} \right)^2\right]^{\frac{1}{2}}(f_{,1})^2\leq
300N\Phi(f_{,1})^2 +$$
$$+
\tfrac{1}{12
N}\sum\left[\left(g-\tfrac{4}{C-f}\right)f_{,i} + \alpha
\tfrac{\rho_{,i}}{\rho}
\right]^{2}(f_{,1})^2+\tfrac{4}{3N}\tfrac{1}{(C-f)^{2}}(f_{,1})^4.
$$
 Note that $K(p^*)=N$. From Step 1 we get
$$ \exp \{N\}\mathcal {A} \leq\tfrac{4}{\alpha}
\exp \left \{-\tfrac{m}{C-f}+ K \right \}\tfrac{g^2
\|\,\text{grad}\; f\|^2}{\rho^{\alpha}}(q^*).$$ It follows that, at
$q^*$,
$$\tfrac{240}{N}\mathcal {A}\exp
\left\{\tfrac{m}{C-f}\right\}\rho^{\alpha}(f_{,1})^2\leq
\tfrac{960}{N\alpha}\;g^{2}(f_{,1})^4.$$
We insert this  estimate into (4.10)
to finally get the third estimate:
$$
4\tfrac{|\rho_{,11}|}{\rho} (f_{,1})^2\leq \tfrac{1}{20}
\tfrac{\|\text{grad}\; K\|^{2}}{K}(f_{,1})^2
+ \tfrac{960}{N\alpha}g^{2}(f_{,1})^4
+ \tfrac{4}{3N}\tfrac{1}{(C-f)^{2}}(f_{,1})^4$$
$$+
302N  \Phi (f_{,1})^2  +
\tfrac{1}{12N}\sum\left[\left(g-\tfrac{4}{C-f}\right)f_{,i} +
\alpha\tfrac{\rho_{,i}}{\rho} \right]^{2}(f_{,1})^2. \eqno(4.11)$$
\vskip 0.1in \noindent After finishing the proof of the three
estimates, we use the inequalities  (4.6), (4.7) and (4.9) and
insert into (4.5); we   get
$$
\tfrac{2}{3N}\sum \left[\left(g - \tfrac{4}{C-f}\right)f_{,i} +
\alpha\tfrac{\rho_{,i}}{\rho} \right]^2(f_{,1})^2 + \left(\alpha -
340N  \right)\Phi (f_{,1})^2$$$$ - 2(N+1) g  (f_{,1})^2  -
2\left(Kg + g -\tfrac{4}{C-f}\right)
\tfrac{\rho_{,1}}{\rho}(f_{,1})^3 $$$$  - \left[(N+1)g' +
\tfrac{1000}{N\alpha  }g^{2}+
\tfrac{4}{3N}\tfrac{1}{(C-f)^2}\right](f_{,1})^4- 328
 \leq 0.\eqno(4.12)$$
As in the proof of Lemma 5.6.11 in \cite{L-X-S-J} we choose $N$ and $\alpha $ such
that
$$1+ N = \tfrac{2\alpha}{3N},\;\;\;i.e., \;\;
\alpha =\tfrac{3N(1+N)}{2},\eqno(4.13)$$ moreover we choose $N$
large enough that $N>10^6,$ and finally we choose\\ $m\geq2C\alpha
N(N+1)$; then
$$g'\,(N+1)\leq \tfrac{1}{N\alpha  }\;g^{2},\qquad \quad
\tfrac{4}{3N}\;\tfrac{1}{(C-f)^{2}}< \tfrac{1}{N^2\alpha}\;g^{2}.$$
 \vskip0.1in \noindent
Again, as in the proof of Lemma 5.6.11 in \cite{L-X-S-J},  we discuss two cases: \vskip
0.1in \noindent {\bf Case 1: \enspace  $\sum \tfrac{\rho_{,i}f_{,i}}{\rho}
> 0$.} \enspace In this case (4.13) gives
the inequality:
$$\tfrac{2}{3N}\left[\left(g - \tfrac{4}{C-f}\right)f_{,i} + \alpha
 \tfrac{\rho_{,i}}{\rho} \right]^2
(f_{,1})^2 - (2 + 2N)\left(g -\tfrac{4}{C-f}\right)
\tfrac{\rho_{,1}}{\rho}(f_{,1})^3
$$$$\geq\tfrac{2}{3N}\left(g
- \tfrac{4}{C-f}\right)^2(f_{,1})^4  \geq
\tfrac{1}{3N}g^2(f_{,1})^4.$$
Note that $K\leq N$, we have
$$2N\left(g - \tfrac{4}{C-f}\right)\tfrac{\rho_{,1}f_{,1}}{\rho}
-2gK\tfrac{\rho_{,1}f_{,1}}{\rho}\geq -
\tfrac{8N}{C-f}\tfrac{\rho_{,1}f_{,1}}{\rho}\geq - 200N\Phi
-\tfrac{3}{50\alpha  }g^{2}(f_{,1})^2.$$
Now we
 insert  the two
estimates  above into (4.12) and get:
$$\tfrac{1}{6N}g^2(f_{,1})^4 - 2(N+1) g  (f_{,1})^2 - 328 \leq 0.$$
\vskip 0.1in \noindent {\bf Case 2: \enspace \enspace $\sum
\tfrac{\rho_{,i}f_{,i}}{\rho} \leq 0$.} \enspace The
inequality of Schwarz implies
$$\tfrac{2}{3N}\sum \left[\left(g -
\tfrac{4}{C-f}\right)f_{,i} +  \alpha  \tfrac{\rho_{,i}}{\rho}
\right]^2(f_{,1})^2 +\tfrac{ \alpha }{4}\Phi(f_{,1})^2 \geq$$
$$\geq
\tfrac{2}{8\alpha +3N}\left(g - \tfrac{4}{C-f}\right)^2
(f_{,1})^4\geq \tfrac{1}{8\alpha +3N}g^2 (f_{,1})^4.$$ We
insert  this estimate into (4.12) and  obtain
the inequality
$$\tfrac{1}{16\alpha +6N}g^2 (f_{,1})^4 - 2(N+1) g  (f_{,1})^2 - 328 \leq 0.$$
Thus, in both cases (1) and (2),
we have an inequality of the type
$$ a_0 g^2 (f_{,1})^4 - a_1g (f_{,1})^2 -
328\leq 0,$$ where \enspace $a_0$ and $a_1$ \enspace are real
positive constants. Consequently
$$\exp \left \{-\tfrac{m}{C-f} + K \right \}\;g^2\tfrac{1}
{\rho^{\alpha} } \; \|\,\text{grad}\; f\|^2\leq a_2,$$ where $a_2$
is a real positive constant depending only on $C$, $\alpha$ and $b$;
this is the upper bound for $\mathcal {D}$
that  we announced in the beginning of Step 2.\\
The upper bound for $\mathcal {D}$ gives an upper bound for
$\mathcal {A}$ (see Step 1), thus  the proof of Proposition 3.8
(Proposition 5.6.13 in \cite{L-X-S-J}) is finished.

In \cite{L-X-S-J}, p. 126, we stated that the proof of  Proposition
5.6.13 there (that is Proposition 3.8 here) is ``similar'' to that
of Proposition 5.6.12 in the same monograph (that is Proposition 3.6
here). We hope that our foregoing proof, containing all details of
the auxiliary tools, will convince the reader
that this statement of ``similarity'' is correct.\\
This finishes this  proof of Chern's Affine Bernstein
Conjecture.
\subsection*{Acknowledgment}
The  authors acknowledge partial support by:\\
{\bf (1)} \enspace the first author by NKBRPC\,(2006CB805905), \enspace NSFC
10631050 \\ and \enspace RFDP\,(20060610004);\\
{\bf (2)} \enspace the second  author by NSFC 10926172;\\
{\bf (3)} \enspace the third   author within the projects DFG PI
158/4-6,  PI 158/8-1 and by four Chinese universities (Peking U,
Tsinghua U, Sichuan U, Yunnan Normal U). Moreover he would like to
thank his colleagues for their great  hospitality during  research
stays in Brazil  (U Brasilia and the  Federal U of  Ceara at
Fortaleza) in summer and at the four Chinese universities in autumn 2010 .\\
The authors acknowledge the interest of X. Wang in the proof of the case $\beta = 0$.

\end{document}